\documentclass[12pt,a4paper]{amsart}

\usepackage{mathrsfs}
\usepackage{amssymb}
\usepackage{amsmath,mathtools}
\usepackage[utf8]{inputenc}
\usepackage{hyperref}

\newtheorem{theorem}{Theorem}[section]
\newtheorem{lemma}[theorem]{Lemma}

\theoremstyle{definition}
\newtheorem{definition}[theorem]{Definition}
\newtheorem{hypothesis}[theorem]{Hypothesis}
\newtheorem{notation}[theorem]{Notation}
\theoremstyle{remark}
\newtheorem{remark}[theorem]{Remark}

\newcommand{\PA}{\mathsf{PA}}
\newcommand{\ZFC}{\mathsf{ZFC}}
\newcommand{\N}{\mathbb N}
\newcommand{\Z}{\mathbb Z}
\newcommand{\D}{\mathcal D}

\newcommand{\Tile}{\operatorname{Tile}}

\title[A true $\Pi^0_1$ sentence  independence of $\PA$]{A true $\Pi^0_1$ sentence  independence of $\PA$ }

\author[M.\  Golshani]{Mohammad Golshani}

\address{Mohammad Golshani, School of Mathematics, Institute for Research in Fundamental Sciences (IPM), P.O.\ Box:
19395--5746, Tehran, Iran.}

\email{golshani.m@gmail.com}

\thanks{The author's research has been supported by a grant from IPM (No. 1405030417). He thanks Erfan Khaniki for many useful discussions on this topic,  for reading an early version of the paper and suggesting the references \cite{krajicek} and \cite{pudlak}.}

\begin{document}

\subjclass[2010]{03F30, 52C23, 03E35}

\keywords{$\Pi^0_1$ sentence, $\PA$, tiling}

\begin{abstract}
We give  a true $\Pi^0_1$ sentence $\varphi$ which is
independent of Peano Arithmetic. Moreover, if $\widehat{\phi}$  denotes its interpretation on
$\omega$ as a sentence of set theory, then $\widehat{\phi}$ is  a  $\Pi^0_1$ sentence independent of $\ZFC$.
\end{abstract}
\maketitle

\section{Introduction}

G\"odel's incompleteness theorems show that there are true arithmetical
sentences which are not provable in Peano arithmetic ($\PA$), however these sentences talk about formal proofs or consistency,
and for years it was open if there are some natural mathematical statements independent of $\PA$.
Paris--Harrington ~\cite{ParisHarrington1977} provided such a  mathematical statement. Their statement
is a variant of finite Ramsey theorem and is  $\Pi^0_2$, and one can show that it is equivalent to the \(1\)-consistency of $\PA$. The question of whether there is a
 natural true $\Pi^0_1$ sentence independent of $\PA$
 was asked by Friedman and Harrington who  emphasized
it as an important problem in logic; see \cite[Section 3]{shelah},
see also \cite{pudlak} and the more recent reference \cite{krajicek}, which explicitly discusses the difficulty of finding natural \(\Pi^0_1\) statements independent of Peano Arithmetic.

In \cite{shelah}, Shelah introduced the method of $(L,n)$-models and  claimed
the existence of a natural $\Pi^0_1$ sentence  independence of $\PA$. Valderrama Hern\'andez~\cite{Valderrama2021}
 found a gap in some of the arguments of Shelah, and  also suggested a way to repair
 them, but it is not clear if his modifications can be formalized in weak systems of arithmetic, hence it still does not fix the above gap.
 Also Harvey
Friedman~\cite{Friedman2014} suggested a solution.

 In~\cite[Dream~2.1]{Shelah2003},
Shelah asks for an analogue of forcing that would establish the
independence of reasonable arithmetical statements from $\PA$. In personal communication \cite{villaveces}, Villaveces suggested that Shelah's method of $(L,n)$-models, or suitable variants of it, may have considerable  potential, particularly in developing forcing techniques for arithmetic. Although the current framework of $(L,n)$-models does not fully achieve this goal, it may provide a promising foundation for further investigation, and we hope to consider this in the future.

In \cite{GT}, Greenfeld and Tao proved a statement independent of $\ZFC$.
Our first interest in this paper was to see if they have obtained the result using forcing or not. After looking at the paper,
it became clear that no forcing argument is used, and indeed the idea is to encode an independent computation into a tiling problem.
After looking at their statement, it became clear that  one can express it as a $\Pi^0_1$ sentence,
and this is the main result of this paper.

\section{a  $\Pi^0_1$ independence sentence}
We fix with the following hypotheses, following \cite{GT}.
For $n \geq 1$ set $\Z_n:=\Z/n\Z$.
\begin{hypothesis}
\begin{enumerate}
\item $G_0=\prod_{i=1}^m\Z_{n_i}$ is a finite abelian group, where $m\in\N$, $n_i\geq1$,

\item $\D=(G_0,E_0,F_1,F_2)$, where $E_0\subseteq G_0$ and
$F_1,F_2\subseteq G:=\Z^2\times G_0$ are finite and nonempty,

\item  $r_{F_1, F_2}\geq1$ is the least integer which bounds $|n|, |m|$, for every
pair $(n, m)$ that appears as integer coordinate of $F_1\cup F_2$. To be more precise,
let
$\pi:G\longrightarrow\Z^2$ be the projection map. Then
$r_{F_1, F_2}=\max\big(\{1\} \cup \{|n|, |m|: \exists f \in F_1 \cup F_2,~ \pi(f)=(n, m)              \}\big)$,

\item  $E:=\Z^2\times E_0$,
\item $ \Tile(\D):=
 \{(A_1,A_2): A_1, A_2 \subseteq G,~(A_1\oplus F_1)\uplus(A_2\oplus F_2)=E\}.$\footnote{In \cite{GT} this is denoted by $\Tile(F_1,F_2;E).$}

\item $\Psi_{\D}$ is the assertion $\Tile(\D)$ is nonempty.
\end{enumerate}
\end{hypothesis}
\begin{notation}
If $\phi$ is a statement in PA, then $\widehat{\phi}$  denotes the interpretation of $\phi$ on
$\omega$ in ZFC.
\end{notation}
\begin{theorem}\label{thm:GT}(see \cite[Theorem~1.8]{GT})
There is  $\D$ as above such that
$\Psi_{\D}$ is independent of ZFC.
\end{theorem}

\begin{proof}
We sketch the proof, following \cite[Sections~3--8]{GT}.
Let $\theta$ be a true $\Pi^0_1$ sentence independent of
$\ZFC$, obtained by the G\"odel--Rosser theorem, and let
$M$ be a Turing machine which halts exactly when $\theta$
fails. By \cite{Berger}, there are finitely many
tiles $T_1,\ldots,T_J\subseteq\Z^2$ such that
there are $A_1,\ldots,A_J\subseteq\Z^2$ such that
$
\Z^2=\biguplus_{j=1}^J(A_j\oplus T_j)
$
 if and only if $M$
does not halt. 

Now the idea is to reduce this tiling problem
to one which involves only two tiles. We briefly recall the
main steps. At the first step, we encode a tiling by $T_1,\ldots,T_J$ by a function
$
c:\Z^2\longrightarrow
\biguplus_{j=1}^J(\{j\}\times T_j)
$
which satisfies
\[
c(x)=(j,t)\Longrightarrow
c(x-t+t')=(j,t'), 
\]
for $t,t'\in T_j.$
Since the range of $c$ is finite, these conditions can
be expressed as a finite system of local Boolean
constraints. By introducing auxiliary functions, these
constraints are reduced to a finite system of linear
functional equations for functions taking values in
$\{-1,1\}$.
The next step is to encode these equations by a finite
system of tiling equations involving two tiles over
$\Z^2\times G'_0$, where $G'_0$ is a finite abelian group.
Finally, by adding a suitable finite cyclic factor to
$G'_0$, all these equations are combined into a single
equation
\[
(A_1\oplus F_1)\uplus(A_2\oplus F_2)
=\Z^2\times E_0.
\]
All the reductions preserve solvability and can be
formalized in $\ZFC$. Thus, for the resulting
$\D=(G_0,E_0,F_1,F_2)$, we have
$
\ZFC\vdash
\Psi_{\D}\leftrightarrow\widehat{\theta}.
$
Since $\widehat{\theta}$ is independent of $\ZFC$,
the same holds for $\Psi_{\D}$.
\end{proof}

\begin{definition}\label{def:patch}
Assume $n\in\N$ and $\D$ is as above. Then
\begin{itemize}
\item $B_n:=[-n,n]^2\cap\Z^2,$
\item  $C_n:=B_n\times E_0,$
\item $D_n:=B_{n+r}\times G_0,$ where $r=r_{F_1, F_2}$,

\item $ \Tile(\D, n)=\{(A_1,A_2): A_1, A_2 \subseteq D_n,~
 C_n\subseteq(A_1\oplus F_1)\uplus(A_2\oplus F_2)\subseteq E\},$

\item $\Psi_{\D, n}$ is the assertion  $\Tile(\D, n)$ is nonempty.
\end{itemize}
\end{definition}
The next lemma follows from compactness theorem.\footnote{In the following proof we use weak K\"onig's  lemma which follows from compactness theorem.}
\begin{lemma}\label{lem:finite-patch-compactness}
Let $\D$ be as in Definition~\ref{def:patch}. Then
\[
\ZFC\vdash \Psi_{\D} \leftrightarrow (\forall n \in \N) \Psi_{\D, n}.
\]
\end{lemma}
\begin{proof}
Work in $\ZFC$.
  By the choice of $r$, if
$x\in C_n$ and $f\in F_1\cup F_2$, then
$x-f\in D_n.$

First assume that $\Psi_{\D}$ holds, and let $A_1,A_2\subseteq G$ be
such that
 $
(A_1\oplus F_1)\uplus(A_2\oplus F_2)=E.
 $
For  $n\in\N$ and $i=1,2$ set $A_i^{(n)}=A_i\cap D_n$.
It is easily seen that
\[
C_n\subseteq
(A_1^{(n)}\oplus F_1)\uplus(A_2^{(n)}\oplus F_2)
\subseteq E.
\]
It follows that $(A_1^{(n)},A_2^{(n)})\in\Tile(\D,n)$, and hence
$\Psi_{\D,n}$ holds for every $n\in\N$.

Conversely, suppose that $\Tile(\D,n)\neq\varnothing$ for every
$n\in\N$. Clearly if $m\geq n$, then
$
\rho_{m,n}:\Tile(\D,m)\longrightarrow\Tile(\D,n),$
defined by
\[
\rho_{m,n}(P_1,P_2)=(P_1\cap D_n,P_2\cap D_n)
\]
is well-defined.
Let $T$ be a rooted tree of height $\omega$, whose $(n+1)$-th level is
\[
\text{Lev}_{n+1}(T)=\{(n,p):p\in\Tile(\D,n)\},
\]
and such that
\[
\text{Suc}_T(n, p)=\{(n+1,q): \rho_{n+1,n}(q)=p     \}.
\]
Since $D_n$ is finite, each level of the tree is finite and by our assumption it is nonempty.
Since each level has a computable size bound,  by weak K\"onig's lemma, it
has an infinite branch
$b= \langle (n, p_n): n \in \omega          \rangle$. Let
$p_n=(A_1^{(n)},A_2^{(n)}).$
Then, for $i=1,2$,
$
A_i^{(n+1)}\cap D_n=A_i^{(n)},
$
 in particular $A_i^{(n)}\subseteq A_i^{(n+1)}$. Set
$
A_i=\bigcup_{n\in\N}A_i^{(n)}.
$
It is evident that $(A_1\oplus F_1)\uplus(A_2\oplus F_2)=E$, hence $(A_1, A_2) \in \Tile(\D)$,
and thus
 $\Psi_{\D}$ holds.
\end{proof}

\begin{lemma}\label{l3}
There is a
$\Delta^0_0$ formula $\Phi_{\D}(n,w)$ in the language of arithmetic such that
if $\varphi_{\D}:= \forall n\,\forall w\,\Phi_{\D}(n,w)$, then
$$\ZFC\vdash
 \widehat\varphi_{\D} \leftrightarrow \Psi_{\D}.$$
 Moreover if $\N\not\models\varphi_{\D}$, then $\PA\vdash\neg\varphi_{\D}$.
\end{lemma}

\begin{proof}
We sketch the proof. Fix $n$ and set
$N=|G_0|(2(n+r)+1)^2=|D_n|$ and $t=2^{2N}$.
There are $t$ many pairs of subsets of $D_n$, and we can
encode them by the $2N$ binary digits of a number
$s<2^{2N}$.
Then $\Psi_{\D,n}$ can be expressed as for some $s<t$,
if $(A_1,A_2)$ is coded by $s$, then the conclusion
\[
C_n\subseteq(A_1\oplus F_1)\uplus(A_2\oplus F_2)
\subseteq E
\]
holds.
Indeed, the $2N$ binary digits of $s$ determine the two
subsets of $D_n$, and the conditions in
Definition~\ref{def:patch} can be checked by bounded
quantification over $D_n$ and the finite sets $F_1,F_2$.
In particular, one checks that every point of $C_n$
is covered, that all representations are unique, and
that no point outside $E$ is covered.

To be more precise, we use G\"odel's $\beta$-function
to provide a code $(b,c)$ for the finite list
$1,2,4,\ldots,t=2^{2N}$.
To check that this is such a code, it suffices to check
that the first entry is $1$, that each successive entry
is twice the preceding one, and that the last entry is
$t$; these are all bounded conditions.
Now let $\Phi_{\D}(n,w)$ be the formula which says that
every such power code with $b,c,t\le w$ admits an $s<t$
satisfying the conditions in Definition~\ref{def:patch}.

To see that $\Phi_{\D}(n,w)$ is $\Delta^0_0$, recall that
the $\beta$-function gives the remainder of $b$ upon
division by $1+(i+1)c$. Since $b,c\le w$, the quotient
and remainder used to determine its $i$-th entry are
both at most $w$. Similarly, the $i$-th binary digit
of $s$ is determined by dividing $s$ by $2^{i+1}$,
as given by the power code, and comparing the remainder
with $2^i$. Since $s<t\le w$, the quotient and remainder
are again at most $w$. Moreover, the indices in the
power code are bounded by $2N$, and the conditions in
Definition~\ref{def:patch} involve only elements of
$D_n$ and $F_1,F_2$. As $F_1,F_2\subseteq D_n$,
these elements can be enumerated using indices less
than $N$. Thus every quantifier in the formula can
be bounded by the term
\[
w+2|G_0|(2(n+r)+1)^2+1,
\]
and hence $\Phi_{\D}(n,w)$ is $\Delta^0_0$.

If $\Psi_{\D,n}$ holds, then any pair in
$\Tile(\D,n)$ gives an $s<t$ for every such power code.
Conversely, if $\Psi_{\D,n}$ fails, take $w$ large
enough to bound $b,c,t$ for a code of the above list.
Then $\Phi_{\D}(n,w)$ fails. Thus
$(\forall w)\Phi_{\D}(n,w)$ is equivalent to
$\Psi_{\D,n}$, and the first conclusion follows from
Lemma~\ref{lem:finite-patch-compactness}.
The rest is clear.
\end{proof}

The following theorem is immediate.
\begin{theorem}\label{thm:main}
Let $\D$ be as in Theorem \ref{thm:GT}. Then the sentence
$\varphi:=\varphi_{\D}$ is true in $\N$, and we have the following:
\begin{enumerate}
\item[(a)] $\varphi$ is independent of $\PA$,
\item[(b)] $\widehat\varphi$ is independent of $\ZFC$.
\end{enumerate}
\end{theorem}

The sentence $\widehat\varphi$ from Theorem \ref{thm:main} is indeed independent of $\ZFC+\mathsf{V=L}.$

\begin{remark}
We can indeed choose $\D$ 
such that $\widehat{\varphi}_{\D}$ is independent
of $\ZFC+\mathsf{V=L}$, and both
$\ZFC+\widehat{\varphi}_{\D}$ and
$\ZFC+\neg\widehat{\varphi}_{\D}$
are equiconsistent with $\ZFC$.

To see this, let $T=\ZFC+\mathsf{V=L}$ and let
$\rho$ be a Rosser sentence for $T$.
By Rosser's  theorem, see
\cite{HP},
both $T+\rho$ and $T+\neg\rho$ are
equiconsistent with $T$, and hence with $\ZFC$.
By the proof of Theorem \ref{thm:GT}, we can choose $\D$ such that
$
\ZFC\vdash\Psi_{\D}\leftrightarrow\widehat{\rho}.
$
The result follows from Lemma~\ref{l3}.
\end{remark}

{\bf AI Declaration.} The author used AI to provide a sketch of the proof of Theorem ~\ref{thm:GT}.

\end{document}